\documentclass[11pt]{article}
\usepackage{amsmath,amsthm,amsfonts,amscd,bezier}
\usepackage{amssymb}
\usepackage[usenames]{color}
\newtheorem{theorem}{Theorem}

\newtheorem{example}[theorem]{Example}
\newtheorem{definition}[theorem]{Definition}
\newtheorem{corollary}[theorem]{Corollary}

\newtheorem{remark}[theorem]{Remark}

\newtheorem{proposition}[theorem]{Proposition}
\newcommand{\pf}{\noindent{\bf Proof:\ \ }}
\newcommand{\n}{\nonumber}
\newcommand{\cqd}{{\hfill $\rule{2mm}{2mm}$}\vspace{1cm}}

\begin{document}

\pagestyle{plain} \pagenumbering{arabic}
\title{Standard Bases for Fractional Ideals of the Local Ring of an Algebroid Curve}
\author{Carvalho, E. and Hernandes, M. E.
\thanks{The first author
was partially supported by CAPES and the second one by CNPq.}\ \
\thanks{Corresponding author: Hernandes, M. E.; email: mehernandes@uem.br }}
\date{ \ }
\maketitle

\begin{center}

key words: Algebroid Curves, Semiring and Semimodules of Values, Standard Bases.
\end{center}

\begin{abstract} In this paper we present an algorithm to compute a Standard Basis for a fractional ideal $\mathcal{I}$ of the local ring $\mathcal{O}$ of an $n$-space
algebroid curve with several branches. This allows us to
determine the semimodule of values of $\mathcal{I}$. When
$\mathcal{I}=\mathcal{O}$, we may obtain a (finite) set of
generators of the semiring of values of the curve, which determines its classical semigroup. In the complex context, identifying the K\"{a}hler differential module
$\Omega_{\mathcal{O}/\mathbb{C}}$ of a plane curve with a
fractional ideal of $\mathcal{O}$ and applying our algorithm, we can
compute the set of values of $\Omega_{\mathcal{O}/\mathbb{C}}$, which
is an important analytic invariant associated to the curve.
\end{abstract}

\section{Introduction}

In \cite{Abramo-Marcelo3} the analytic classification problem for irreducible plane curves was solved using discrete structures, namely the semigroup of values and the set of values of K\"{a}hler differentials of the local ring of the curve.

The semigroup of an irreducible analytic plane curve (also called a
branch) is a classical object. It is finitely generated and its
minimal set of generators can be computed via several methods: by
characteristic exponents of a Puiseux expansion, by the sequence of
multiplicity of the infinitely near points of the canonical
resolution of the branch, by intersection multiplicities of the
branch with maximal contact curves or with approximate roots, etc.
(See \cite{hefez}, \cite{popescu} and \cite{zariski}.)

Natural generalizations consist in considering algebroid curves, space irreducible curves or plane curves with several branches.

The authors of \cite{castellanos1}, \cite{castellanos2} and \cite{Abramo-Marcelo} present distinct methods to obtain the generators of the semigroup of values of an irreducible algebroid space curve: by the infinitely near points associated to any resolution (see \cite{castellanos1} and \cite{castellanos2}) and by Standard Bases for the local ring of the curve (see \cite{Abramo-Marcelo}).

The semigroup of values $S$ of an algebroid plane curve $Q$ with $r$ branches is a (not finitely generated) submonoid of $\mathbb{N}^r$. Delgado in \cite{delgado} shows that $S$ is determined by a finite subset of (relative maximal) points of $S$ and all the semigroups of curves with $r-1$ branches obtained from $Q$.

Considering an abstract (good) semigroup $S$ of $\mathbb{N}^r$, so $S$ is not necessarily associated to a curve (see \cite{barucci2}), D'Anna {\it et al.} (see \cite{DGMT}) present another finite set of elements that determines $S$.

In this paper we take into account all the situations above, that is, we consider an $n$-space algebroid curve $Q$ with $r$ branches.

The elements of the semigroup $S$ of $Q$ correspond to values of non zero divisors of the local ring $\mathcal{O}$ of $Q$. In \cite{carvalho}, we consider the set of values $\Gamma$ of {\it all} elements in $\mathcal{O}$ and we prove that $\Gamma$ is a finitely generated semiring of $\overline{\mathbb{N}}^r:=(\mathbb{N}\cup\{\infty\})^r$ equipped with the tropical operations. Maugeri and Zito in \cite{zito} consider this structure of semiring to study the embedding dimension of a good semigroup.

For an analytic plane curve $Q$ with $r$ branches, Zariski in
\cite{equising} characterizes the topological class (the
equisingularity type) of $Q$ as embedded curve by means of the semigroup of each branch
and the mutual intersection multiplicity of the branches. On the
other hand, Waldi in \cite{waldi} shows that the semigroup of $Q$
determines completely its topological class. Using the semiring
structure of $\Gamma$, its minimal set of generators is
given by a finite set $A$ of points in $\mathbb{N}^r$, which Delgado
in \cite{delgado} calls irreducible absolute maximal points, and a
set $B$ with $r$ points that correspond to the intersection
multiplicity of a branch with the others. As $A$ and $B$ can be
obtained from the semigroup of each branch and the mutual
intersection multiplicity of the branches, this approach allows us
to directly connect the results of Zariski and Waldi, because
$S=\Gamma\cap\mathbb{N}^r$.

In the study of analytic classification of a plane curve $Q$ another algebraic structure has important role: the K\"{a}hler differential module
$\Omega_{\mathcal{O}/\mathbb{C}}$ of the local ring $\mathcal{O}$.
We can consider $\Omega_{\mathcal{O}/\mathbb{C}}$ as a fractional ideal of
$\mathcal{O}$ and its set of values $\Lambda$ allows us to stratify the topological class. This is an approach that has been used in the analytic classification
problem (see \cite{Abramo-Marcelo3} for the irreducible case and
\cite{AME} for curves with two branches).

Several authors have considered the set of values of fractional ideals of $\mathcal{O}$ and they obtained interesting results related with symmetry properties and their codimensions: D'Anna (in \cite{danna}); Barucci, D'Anna and Fr\"oberg (in \cite{barucci}); Guzm\'an and Hefez (in \cite{edison}); Pol (in \cite{pol2}), etc.

Pol in \cite{pol} investigated the set $v(\mathcal{R}_Q)$ of values of the module $\mathcal{R}_Q$ of logarithmic residues along a complete intersection curve $Q$. In this case, $\mathcal{R}_Q$ coincides with the dual of the Jacobian ideal $J(Q)$ of $Q$. In particular, we have that $v(\mathcal{R}_Q)$ and the set $v(J(Q))$ of values of $J(Q)$ determine each other and $v(J(Q))=\sigma -\underline{1}+\Lambda$, where $\sigma$ is the conductor of the semigroup $S$ (see Theorem 1.2 and Proposition 3.31, \cite{pol}). In addition, in Subsection 4.3.3 of \cite{pol} it is presented an algorithm to compute the set $v(\mathcal{R}_Q)$ or, equivalently, $v(J(Q))$ and $\Lambda$ for plane curves with two branches.

The aim of this work is to develop a Standard Bases theory for
fractional ideals $\mathcal{I}$ of the local ring of an $n$-space
algebroid curve $Q$ with $r$ branches and to present an algorithm to
compute a Standard Basis for $\mathcal{I}$ and, consequently, a
finite set of generators for the semimodule of values of
$\mathcal{I}$.

This paper is organized as follows. In Section 2 we present the
concepts of semiring and semimodule associated to the local ring and
its fractional ideals. We develop the theory of Standard Bases for
a fractional ideal in Section 3 in order to obtain a finite set of
generators for the semimodule associated. The K\"{a}hler
differential module $\Omega_{\mathcal{O}/\mathbb{C}}$ of the local
ring of a plane curve is considered in Section 4 to illustrate an application
of the algorithm. In particular, we are able to compute the set of
values of $\Omega_{\mathcal{O}/\mathbb{C}}$ for an algebroid curve
with several branches without restrictions on the number of branches in contrast to the algorithm presented in
\cite{pol} which works only for curves with two irreducible
components.

\section{Fractional and Relative Ideals}

Let $\mathbb{K}$ be an arbitrary algebraically closed field and let $(\mathcal{O},\mathcal{M},\mathbb{K})$ be a one dimensional complete Noetherian local ring containing $\mathbb{K}$ as a coefficient field. If $\mathcal{O}$ has embedding dimension $n$, that is, $n=dim_{\mathbb{K}}\frac{\mathcal{M}}{\mathcal{M}^2}$, then $\mathcal{O}\simeq\frac{\mathbb{K}[[X_1,\ldots ,X_n]]}{Q}$ for some proper radical ideal $Q\subset \mathbb{K}[[X_1,\ldots ,X_n]]$.

In this paper, by an algebroid curve in $\mathbb{K}^n$, or an $n$-space curve, we mean $Spec(\mathcal{O})$ with $\mathcal{O}=\frac{\mathbb{K}[[X_1,\ldots ,X_n]]}{Q}$ as above. In what follows we will simply say ``a curve $Q$'' and we will call $\mathcal{O}$ the local ring of the curve.

If $Q=\bigcap_{i=1}^{r}P_i$, where $P_i\subset \mathbb{K}[[X_1,\ldots ,X_n]]$ is a prime ideal, then we call $P_i$ a {\it branch} of the curve $Q$. The integral closure $\overline{\mathcal{O}_i}$ of the local ring $\mathcal{O}_i=\frac{\mathbb{K}[[X_1,\ldots ,X_n]]}{P_i}$ of each branch in its quotient field $\mathcal{Q}_i$ is a discrete valuation ring.
Denoting $v_i:\mathcal{Q}_i\rightarrow\overline{\mathbb{Z}}:=\mathbb{Z}\cup\{\infty\}$
the normalized discrete valuation associated to
$\overline{\mathcal{O}_i}$, where $v_i(0)=\infty$, and fixing a uniformizing parameter $t_i\in\overline{\mathcal{O}_i}$ for $v_i$, we can consider $\overline{\mathcal{O}_i}=\mathbb{K}[[t_i]]$ and $v_i(h_i)=ord_{t_i}(h_i)$, with $h_i\in\overline{\mathcal{Q}_i}$.

In the sequel we denote $I=\{1,\ldots ,r\}$.

If $\overline{\mathcal{O}}$ is the integral closure of $\mathcal{O}$
in its total ring of fractions $\mathcal{Q}$, then we have
\begin{align}\label{inclusions} \mathcal{O} \hookrightarrow
\bigoplus_{i\in I}\mathcal{O}_i \hookrightarrow
\overline{\mathcal{O}} \simeq\bigoplus_{i\in
I}\overline{\mathcal{O}_i} = \bigoplus_{i\in I}\mathbb{K}[[t_i]]
\hookrightarrow \bigoplus_{i\in I}\mathbb{K}((t_i)) \simeq
\bigoplus_{i\in I}\mathcal{Q}_i = \mathcal{Q}, \end{align} where all
the homomorphisms are injective and the first one is given by
$g\mapsto (g_1,\ldots ,g_r)$ in such way that $g_i$ denotes the
canonical image of $g\in\mathcal{O}$ in $\mathcal{O}_i$.

In what follows we
identify $X_j+P_i\in\mathcal{O}_i$ with its isomorphic image
$x_j(t_i)\in\mathbb{K}[[t_i]]$ and $\mathcal{O}_i$ with
$\mathbb{K}[[x_1(t_i),\ldots ,x_n(t_i)]]$. We call $(x_1(t_i),\ldots
,x_n(t_i))$ a parameterization of $P_i$. The set $$\Gamma_i= v_i(\mathcal{O}_i):=\{v_i(g_i);\
g_i\in\mathcal{O}_i\}\subseteq \overline{\mathbb{N}}$$ is an additive monoid setting $\gamma_i+\infty=\infty$ for all $\gamma_i\in\Gamma_i$. Notice that $S_i=\Gamma_i\cap\mathbb{N}$ is the classical semigroup of values associated to the branch $P_i$.

Let $\psi:\mathbb{K}[[X_1,\ldots ,X_n]]\rightarrow\bigoplus_{i\in
I}\mathbb{K}[[t_i]]=\overline{\mathcal{O}}$ be the homomorphism ring
defined by
$$X_i\mapsto (x_{i}(t_1),\ldots ,x_{i}(t_r)).$$

As $ker(\psi)=Q$, we have
$$\mathcal{O}=\frac{\mathbb{K}[[X_1,\ldots ,X_n]]}{Q}\simeq
Im(\psi) \subseteq \bigoplus_{i\in I}\mathbb{K}[[x_1(t_i),\ldots
,x_n(t_i)]]\simeq \bigoplus_{i\in I} \mathcal{O}_i,
$$ where the induced isomorphism is
given by
 $$g\in\mathcal{O}\ \mapsto\ (g(x_1(t_1),\ldots ,x_n(t_1)),\ldots ,g(x_1(t_r),\ldots ,x_n(t_r))).$$

Identifying $\mathcal{O}$ with a subalgebra of $\bigoplus_{i\in
I}\mathbb{K}[[t_i]]\subset\bigoplus_{i\in I}\mathcal{Q}_i$ and
considering the canonical projection $\pi_J:\bigoplus_{i\in
I}\mathcal{Q}_i\rightarrow\bigoplus_{j\in J}\mathcal{Q}_j$, for a
non empty subset $J$ of $I$, we have that
$\pi_J(\mathcal{O})\approx\frac{\mathbb{K}[[X_1,\ldots
,X_n]]}{\cap_{j\in J}P_j}=:\mathcal{O}_J$ and we have
$\mathcal{O}_{\{i\}}=\mathcal{O}_i$.

By (\ref{inclusions}) we have
$$\Gamma=v(\mathcal{O}):=\{v(g):=(v_1(g),\dots ,v_r(g));\
g\in\mathcal{O}\}\subseteq\bigoplus_{i\in I}\Gamma_i.$$

In particular, the semigroup $S=\Gamma\cap\mathbb{N}^r$ of $Q$ and the set $\Gamma$ determine each other.

For $\alpha = (\alpha_1, \ldots , \alpha_r),\beta = (\beta_1, \ldots
, \beta_r)\in\Gamma$ the following properties are immediate:
\begin{itemize}
\item[{\bf P.0)}] If $\alpha_i=0$ for some $i\in I$,
then $\alpha=\underline{0}:=(0,\ldots ,0)$.
\item[{\bf P.1)}] If $\alpha_k=\beta_k < \infty$ for
some $k\in I$, then there exists $\gamma = (\gamma_1, \ldots ,
\gamma_r) \in\Gamma$ such that
$\gamma_i\geq\min\{\alpha_i,\beta_i\}$ for all $i\in I$ (the
equality holds if $\alpha_i\neq\beta_i$) and
$\gamma_k>\alpha_k=\beta_k$.
\item[{\bf P.2)}] $\inf\{\alpha,\beta\}:=(\min\{\alpha_{1},\beta_{1}\},\ldots
,\min\{\alpha_{r},\beta_{r}\})\in\Gamma$.
\end{itemize}

Moreover, $\Gamma$ is a commutative semiring equipped with the
tropical operations
$$\alpha \oplus\beta=\inf\{\alpha,\beta\}\ \ \ \ \ \text{and} \
\ \ \ \  \alpha\odot\beta=\alpha+\beta.$$

Remark that $(\Gamma, \oplus)$ is an idempotent monoid with element identity $\underline{\infty}:=(\infty, ..., \infty)$.

\begin{definition} We call $(\Gamma,\oplus,\odot)$ the {\rm semiring of
values} associated to the curve $Q = \bigcap_{i\in I}P_i$.
\end{definition}

In \cite{carvalho} we proved that $\Gamma$ is a finitely generated
semiring, i.e., there exists a finite subset $\{\gamma_1,\ldots
,\gamma_m\}\subset\Gamma$ that allows us to express any $\gamma\in\Gamma$ as
\begin{equation}\label{elemento}\gamma=
\left (\gamma_1^{\alpha_{11}}\odot\ldots\odot \gamma_m^{\alpha_{1m}}\right )
\oplus\ldots \oplus \left (\gamma_1^{\alpha_{s1}}\odot\ldots \odot
\gamma_m^{\alpha_{sm}}\right )=\inf\left \{
\sum_{j=1}^{m}\alpha_{1j}\gamma_j,\ldots
,\sum_{j=1}^{m}\alpha_{sj}\gamma_j\right\},\end{equation} with $\alpha_{ij}\in\mathbb{N},\  1\leq j\leq m$ and $1\leq i\leq s$, for some $s\leq r$ which depends on $\gamma$.

Since $\overline{\mathcal{O}}$ is an $\mathcal{O}$-module of finite
type, its conductor ideal
$\mathcal{C}=(\mathcal{O}:\overline{\mathcal{O}})$ (that is an ideal of $\mathcal{O}$ as well) contains a nonzero
divisor and $\mathcal{C}=(t_1^{\sigma_1},\ldots
,t_r^{\sigma_r})\overline{\mathcal{O}}$. We call
$\sigma=(\sigma_1,\ldots ,\sigma_r)\in S\subset\Gamma$ the {\it
conductor} of $\Gamma$.

\begin{remark}\label{condutor} D'Anna (see \cite{danna}) shows that we can determine $\sigma_i$ for all $i\in I$ via the prime ideals $\{P_1,\ldots ,P_r\}$. More explicitly, denote $Q^{i}$ the canonical image of the ideal
$\bigcap_{j\in I\atop j\neq i}P_j+P_i$ of $\mathbb{K}[[X_1, ...,
X_n]]$ in $\mathcal{O}_i$ and $v_i(Q^i)=\{v_{i}(q);\ q\in Q^i\}$. As
$v_i(Q^i)$ is a $\Gamma_i$-monomodule, that is,
$\Gamma_i+v_i(Q^i)\subseteq v_i(Q^i)$, there exists $\delta_i\in
v_i(Q^{i})$ such that $\delta_i+\overline{\mathbb{N}}\subseteq
v_i(Q^{i})$ and $\delta_i-1\not\in v_i(Q^{i})$. By Proposition 1.3 in
\cite{danna}, we have that
$\sigma_i=\delta_i$ for all $i\in I$.
\end{remark}

We can obtain $\sigma_i$ by the algorithm presented in \cite{Abramo-Marcelo} that allows to compute the set of values of any finitely generated $\mathcal{O}_i$-module. If $n=2$, that is, for a plane curve $Q=\langle \prod_{i\in I}f_i\rangle=\cap_{i\in I}\langle f_i\rangle$,
we have that $v_i(Q^{i})=v_i\left( \prod_{j\in I\atop j\neq i}f_j\right)+\Gamma_i
= \sum_{j\in I\atop j\neq i}v_i(f_j)+\Gamma_i$ and
$\sigma_i=\sum_{j\in I\atop j\neq i}v_i(f_j)+\mu_i$, where $\mu_i$ is
the conductor of $\Gamma_i$ that can be computed in terms of the
minimal set of generators of $\Gamma_i$ (or $S_i$).

In a more general situation, we will consider a finitely generated (regular) fractional ideal $\mathcal{I}$ of $\mathcal{O}$.

As before, considering the canonical projection
$\pi_J:\bigoplus_{i\in I}\mathcal{Q}_i\rightarrow\bigoplus_{j\in
J}\mathcal{Q}_j$ with $\emptyset\neq J\subset I$ we have that
$\mathcal{I}_J:=\pi_J(\mathcal{I})$ is a fractional ideal of
$\mathcal{O}_J$ and we denote $\mathcal{I}_{\{i\}}$ by
$\mathcal{I}_i$.

In a natural way, we define the set of values of $\mathcal{I}_i$
(respectively $\mathcal{I}$) by
$$v_i(\mathcal{I}_i)\!:=\!\{v_i(h_i);\
h_i\in\mathcal{I}_i\}\! \subseteq \overline{\mathbb{Z}}\ \ \left
(\mbox{resp.}\  v(\mathcal{I})\!:=\!\{v(h)\!:=(v_i(h_i),\ldots
,v_r(h_r));\ h \in  \mathcal{I}\}\!\subseteq \bigoplus_{i\in I}
v_i(\mathcal{I}_i)\right )\!.$$

The properties {\bf P.1} and {\bf P.2} remain true for elements in $v(\mathcal{I})$ and we also have the following:
\begin{itemize}
    \item[{\bf P.3})] $\Gamma + v(\mathcal{I}) \subseteq v(\mathcal{I})$.
    \item[{\bf P.4})] There exists $\gamma \in \Gamma$ such that $\gamma + v(\mathcal{I}) \subseteq \Gamma$.
\end{itemize}

The properties {\bf P.2} and {\bf P.3} indicate that
$v(\mathcal{I})$ admits a ``conductor'' $\varrho \in v(\mathcal{I})$, that is, $\varrho +
\overline{\mathbb{N}}^r \subseteq v(\mathcal{I})$ and $\varrho
-e_i\not\in v(\mathcal{I})$ for all $e_i$ in the canonical
$\mathbb{Q}$-basis of $\mathbb{Q}^r$. In particular, if $\mathcal{I}
= \mathcal{O}$, then $v(\mathcal{I} )=v(\mathcal{O})=\Gamma$ and
$\varrho =\sigma$.

\begin{definition} The set $v(\mathcal{I})$ (resp. $v(\mathcal{I}_i)$) is called the {\rm relative ideal} associated to $\mathcal{I}$ (resp. $\mathcal{I}_i$).
\end{definition}

In the same way done for $\Gamma$, properties {\bf P.2} and {\bf
P.3} allow us to consider tropical operations in $v(\mathcal{I})$
defined by $$\lambda \oplus\lambda' :=\inf\{\lambda,\lambda'\} \ \ \
\ \ \text{and}\ \ \ \ \ \gamma \odot \lambda := \gamma + \lambda,\ \
\ \ \ \mbox{for all}\ \lambda,\lambda'\in v(\mathcal{I})\
\mbox{and}\ \gamma\in\Gamma.$$

It is immediate that $(v(\mathcal{I}),\oplus)$ is a commutative monoid and, for all $\gamma, \gamma' \in \Gamma$ and $\lambda, \lambda' \in v(\mathcal{I})$, we have the following:
\begin{enumerate}
\item $\gamma \odot (\lambda \oplus \lambda') = (\gamma \odot \lambda) \oplus (\gamma \odot \lambda')$;
\item $(\gamma \oplus \gamma') \odot \lambda = (\gamma \odot \lambda) \oplus (\gamma' \odot \lambda)$;
\item $(\gamma \odot \gamma') \odot \lambda = \gamma \odot (\gamma' \odot \lambda)$;
\item $\underline{0} \odot \lambda = \lambda$.
\end{enumerate}
Hence, these above properties give a $\Gamma$-$\ \!$semimodule structure to $v(\mathcal{I})$.

In the next section we will see that $v(\mathcal{I})$ is a finitely generated $\Gamma$-$\ \!$semimodule, that is, there exists a subset $\{\lambda_1, ..., \lambda_s\} \subset v(\mathcal{I})$ such that any $\lambda \in v(\mathcal{I})$ can be expressed as $$\lambda = (\gamma_1 \odot \lambda_1) \oplus \ldots \oplus (\gamma_s \odot \lambda_s), \ \text{with} \ \gamma_j \in \Gamma \ \text{for all} \ j = 1, ..., s.$$

We will also provide an algorithm to compute a finite subset $H$ of $\mathcal{I}$ such that $v(H):=\{v(h); \ h \in H
\}$ is a set of generators of $v(\mathcal{I})$ as $\Gamma$-$\
\!$semimodule. In particular, the algorithm allows us to compute a finite
set of generators for the semiring $\Gamma$ when we apply it for
$\mathcal{I}=\mathcal{O}$ whose existence was proved in
\cite{carvalho}.

\section{Standard Bases for a Fractional Ideal}

Let $\lambda = (\lambda_1, \ldots , \lambda_r)$ and $\lambda' = (\lambda'_1, \ldots , \lambda'_r)$ be elements of $\overline{\mathbb{Z}}^r$. We consider in $\overline{\mathbb{Z}}^r$ the following product order: $\lambda \leq \lambda'$ if and only if $\lambda_i \leq \lambda'_i$, for all $i \in I = \{1, ..., r\}$. We write $\lambda < \lambda'$ if $\lambda \leq \lambda'$ and $\lambda \neq \lambda'$.

In what follows we will extend the main concepts introduced in
\cite{Abramo-Marcelo} and in \cite{carvalho} for a fractional ideal $\mathcal{I}$ of $\mathcal{O}$.

\begin{definition}\label{reduction}
Let $G\subset \mathcal{M}\setminus\{0\}\subset\mathcal{O}$ and
$H\subset\mathcal{I}\setminus\{0\}$ be finite subsets. Given a
nonzero element $f\in\mathcal{Q}$ and $k \in I_f:=\{i \in I; \
v_i(f)\neq \infty \}$, we say that $f'$ is a $k$-{\rm reduction} of
$f$ modulo $(H, G)$ if there exist $c \in \mathbb{K}$, $h \in H$ and
a $G$-product $G^{\alpha}$, that is, $G^{\alpha}=\prod_{g_i\in
G}g_i^{\alpha_i}$ with $\alpha_i\in\mathbb{N}$, such that $$ f' = f
- cG^{\alpha}h,$$ with $v_i(f') \geq v_i(f)$ for all $i \in I$ and
$v_k(f') > v_k(f)$. We say that  $f'$ is a {\rm reduction} of $f$
modulo $(H,G)$ if $f'$ is a $k$-reduction of $f$ modulo $(H,G)$ for
some $k \in I_f$.
\end{definition}

\begin{remark} Notice that $f\in\mathcal{Q}$ has a reduction modulo $(H, G)$ if and
only if for some $k \in I_f$ there exist a $G$-product $G^{\alpha}$
and $h \in H$ such that $v_i(f) \leq
v_i(G^{\alpha}h)$ for all $i \in I$ and $v_k(f)=v_k(G^{\alpha}h)$.
\label{corolarioReducaoSB}
\end{remark}

Now we introduce the main notion of this paper.

\begin{definition}\label{SB}
A \emph{Standard Basis} for $\mathcal{I}$ is a pair of finite
subsets $G\subset \mathcal{M}\setminus\{0\}$ and
$H\subset\mathcal{I}\setminus\{0\}$ such that every nonzero element
in $\mathcal{I}$ admits a reduction modulo $(H,G)$.
\end{definition}

In \cite{carvalho} we presented the concepts of reduction and
Standard Basis (and its existence as well) for the local ring
$\mathcal{O}$ of an algebroid curve. We can recover these concepts
setting $H = \{1\}$ in Definition \ref{reduction} and Definition
\ref{SB}. In this case, if $(\{1\},G)$ is a Standard Basis for
$\mathcal{O}$, then we will indicate it by $G$ to agree with
the definition presented in \cite{carvalho}.

In order to simplify the next results and their proofs, we will always consider $G$ as a Standard Basis for $\mathcal{O}$. In addition, as $G$ is fixed, we say a Standard Basis $H$ for a fractional ideal $\mathcal{I}$ instead $(H,G)$.

The following result guarantees the existence of a Standard Basis
$H$ for any fractional ideal $\mathcal{I}$ of $\mathcal{O}$.
Although its proof is quite similar to the proof of the existence of
a Standard Basis for $\mathcal{O}$ given in \cite{carvalho}, we will present it to introduce
certain sets that will be important in the algorithm that we will
develop to compute $H$.

\begin{theorem}\label{existBase2}
    Any fractional ideal $\mathcal{I}$ of the local ring $\mathcal{O}$ of an
    algebroid curve $Q = \bigcap_{i\in I}P_i$ admits a Standard Basis.
\end{theorem}

\pf For each $i \in I$, we denote by $\mathcal{I}^i$ the canonical
image of the $\mathcal{O}$-module $$\mathcal{T}^i = \{f \in
\mathcal{I}; \ v_j(f) = \infty, \ \text{for all}\ j \in I \setminus
\{i\} \}$$ in $\mathcal{I}_i$ and we consider $B_i \subset
\mathcal{I}$ such that the homomorphic image of $B_i$ in
$\mathcal{I}_i$ is a Standard Basis for $\mathcal{I}^i$, which can
be computed as described in \cite{Abramo-Marcelo}. As the
homomorphic image of any finite subset $A$ of $\mathcal{I}$ such
that $v_i(A) = v_i(B_i)$ is a Standard Basis for $\mathcal{I}^i$, we
can take $B_i$ as a subset of $\mathcal{T}^i$, that is,
$v_j(h)=\infty$ for all $h \in B_i$ and $j\in I\setminus\{i\}$.

Let $B_0'$ be a subset of $\mathcal{I}$ satisfying $v(B_0') =
v(\mathcal{I})$ such that $v(h) \not\in v(B_0'\backslash \{h\})$ for
all $h \in B_0'$ and set $B_0:=\{h \in B_0'\setminus\{0\};\ v_i(h)
<\kappa_i \ \mbox{if} \ i\in I_h\}$, where $\kappa_i$ is the
conductor of $v_i(\mathcal{I}^i)$.

Now consider the finite set $H=\bigcup_{i=0}^{r}B_i$. We will show
that $H$ is a Standard Basis for $\mathcal{I}$.

Let $f$ be a nonzero element in $\mathcal{I}$. If $v_i(f)<\kappa_i$ for all $i \in I_f$, then there exists $h \in B_0$ such that $v(f)= v(h)$. If $\kappa_k\leq v_k(f)$ for some $k \in I_f$, then $v_k(f)\in v_k(\mathcal{I}^k)$. As the homomorphic images of $G \subset \mathcal{O}$ and $B_k\subset\mathcal{I}$ are Standard Bases for
$\mathcal{O}_k$ and $\mathcal{I}^k$ respectively, there exist a $G$-product $G^{\alpha}$ and $h_k \in B_k$ such that $v_k(f)=v_k(G^{\alpha}h_k)$ and $v_i(f)\leq v_i(G^{\alpha}h_k)=\infty$ for all $i\in I\setminus\{k\}$.

By Remark \ref{corolarioReducaoSB}, we conclude that $H$ is a
Standard Basis for $\mathcal{I}$. \cqd

Notice that in the proof of the above theorem we can change
$(\kappa_1,\ldots ,\kappa_r)=:\kappa$ for any upper bound of itself.

\begin{remark}\label{Conjuntos B_i's}
For $\mathcal{I}=\mathcal{O}$ we have $\mathcal{I}^i=Q^i$. In this
way, if we consider $B_i$ as the union of Standard Bases for
$\mathcal{O}_i$ and $Q^i$ for all $i\in I$, then the same proof of
the above theorem shows that $\bigcup_{i=0}^{r}B_i$ is a Standard
Basis for $\mathcal{O}$. It is important to emphasize that $B_i$, for
$i\in I$, can be computed as described in \cite{Abramo-Marcelo}.
\end{remark}

Theorem \ref{existBase2} allows us to conclude that the relative ideal associated to a fractional ideal $\mathcal{I}$ of $\mathcal{O}$ is finitely generated as a $\Gamma$-semimodule.

\begin{corollary}\label{finitgerad}
    The relative ideal $v(\mathcal{I})$ associated to the fractional ideal
    $\mathcal{I}$ is a $\Gamma$-semimodule generated by $v(H)$,
    where $H$ is a Standard Basis for $\mathcal{I}$.
\end{corollary}
\pf It is immediate of the previous theorem and Remark
\ref{corolarioReducaoSB}.
\cqd

In the sequel we will start to develop the necessary tools and
results in order to obtain an algorithm to compute a Standard Basis for a
fractional ideal $\mathcal{I}$ of $\mathcal{O}$. Consequently, this
algorithm will provide a finite set of generators for the relative ideal
$v(\mathcal{I})$.

Let $H$ be a subset of $\mathcal{I}\setminus\{0\}$. Given
$f\in\mathcal{Q}\setminus\{0\}$, we have a sequence (possibly
infinite) of reductions modulo $(H,G)$: $$f_0=f,\ \ \ \
f_i=f_{i-1}-c_iG^{\alpha_i}h_i,\ i > 0,$$ where $c_i\in\mathbb{K}$,
$G^{\alpha_i}$ is a $G$-product and $h_i \in H$.

We say that $f'$ is a {\it final reduction} of $f \in \mathcal{Q}\setminus
\{0\}$ modulo $(H,G)$ if $f'$ is obtained from $f$ via a sequence
of reductions modulo $(H,G)$ and there does not exist a reduction of
$f'$ modulo $(H,G)$.

In case of infinite reductions, we get a sequence
$s_k=\sum_{i=1}^{k}c_iG^{\alpha_i}h_i,\ k\geq 1$, in $\mathcal{I}$.
As we have $v(G^{\alpha_i}h_i)\neq v(G^{\alpha_j}h_j)$ for $i\neq
j$, the set $\{c_iG^{\alpha_i}h_i;\ i\geq 1\}$ is summable and the
sequence $s_k$ is convergent in $\mathcal{I}$.

If $H$ is a Standard Basis for $\mathcal{I}$, then every element
$f\in\mathcal{I}\setminus\{0\}$ admits a sequence of reductions
modulo $(H,G)$ to $0$, that is,
$f=\lim_{k\rightarrow\infty}\sum_{i=1}^{k}c_iG^{\alpha_i}h_i$ or,
equivalently, we can write $$f = \sum_{\delta \in
\Delta}c_{\delta}G^{\delta}h_{\delta},$$ where $\Delta$ is a subset
of $\mathbb{N}^{\sharp G}$, $c_{\delta} \in \mathbb{K}$,
$G^{\delta}$ is a $G$-product, $h_{\delta} \in H$ and $v(f) \leq
v(G^{\delta}h_{\delta})$, for all $\delta \in \Delta$. In this case,
$H$ is a set of generators for $\mathcal{I}$. Moreover, in order to
verify if an element of $\mathcal{Q}$ belongs to $\mathcal{I}$ it is
sufficient that it has a vanishing final reduction modulo $(H,G)$.

\begin{remark}\label{Obs_Altura}
Let $H\!\subset \mathcal{I}$ and $f=\sum_{\delta \in
\Delta}c_{\delta}G^{\delta}h_{\delta}\in\mathcal{I}$ with
$h_{\delta}\!\in\! H$.\! If\ \ $\inf\{v(G^{\delta}h_{\delta}); \delta \in
\Delta \} = v(f)$, then for each $k\in I$ there exists $\delta_{k}
\in \Delta$ such that $v_k(G^{\delta_{k}}h_{\delta_{k}}) =
\min\{v_k(G^{\delta}h_{\delta}); \ \delta \in \Delta \} = v_k(f)$.
For $i \in I\setminus \{k\}$, we have
$$v_i(f) = \min\{v_i(G^{\delta}h_{\delta}); \ \delta \in
\Delta \} \leq v_i(G^{\delta_{k}}h_{\delta_{k}}).$$
\end{remark}

It is obvious that for any Standard Basis $H$ for $\mathcal{I}$ and
for every $h\in\mathcal{I}\setminus\{0\}$, the set $H \cup \{h\}$ is
also a Standard Basis for $\mathcal{I}$. So it is natural to
introduce the following definition.

\begin{definition}
    Let $H$ be a Standard Basis for $\mathcal{I}$. We say that $H$ is
    {\rm minimal} if for every $h \in H$ there does not exist a reduction of
    $h$ modulo $(H\setminus\{h\},G)$.
\end{definition}

The guarantee of existence of a minimal Standard Basis is proved in
the sequel. More precisely, we show that from a Standard Basis $H$
we can obtain a minimal Standard Basis discarding elements $h\in H$
that admit some reduction modulo $(H\setminus\{h\},G)$.

\begin{proposition}\label{minima}
    Let $H$ be a Standard Basis for $\mathcal{I}$. If $h\in H$ admits some reduction
    modulo $(H\setminus\{h\},G)$, then $H' = H\setminus\{h\}$ is a Standard Basis for $\mathcal{I}$.
\end{proposition}

\pf Suppose that $h\in H$ admits a $k$-reduction modulo $(H',G)$ for
some $k \in I_h$, that is, there exist $c_1 \in \mathbb{K}$, a
$G$-product $G^{\alpha_1}$ and $h_1 \in H'$ such that $h' = h -
c_1G^{\alpha_1}h_1$ satisfies $v_i(h') \geq v_i(h)$ for all $i \in
I$ and $v_k(h') > v_k(h)$. In particular,
$v_k(h)=v_k(G^{\alpha_1}h_1)$.

Given $f\in\mathcal{I}\setminus\{0\}$, suppose that
$f'=f-cG^{\alpha}h$ is a $j$-reduction of $f$ modulo $(H,G)$, that
is, $v_i(f')\geq v_i(f)$ for all $i\in I$ and $v_j(f')>v_j(f)$.

If $v_j(h) = v_j(G^{\alpha_1}h_1)$, then there exists $c'\in\mathbb{K}$ such that $f'=f-c'G^{\alpha+\alpha_1}h_1$ is a $j$-reduction of $f$ modulo $(H',G)$.

On the other hand, that is, if $v_j(h) < v_j(G^{\alpha_1}h_1)$, then
we have $v_j(h') = v_j(h)$ and, consequently, $j\in I_{h'}$, which
implies that $h'$ admits a $j$-reduction modulo $(H,G)$. In this
way, there exist $c_2 \in \mathbb{K}$, a $G$-product $G^{\alpha_2}$
and $h_2 \in H$ such that $h'' = h' - c_2G^{\alpha_2}h_2 = h
-c_1G^{\alpha_1}h_1 - c_2G^{\alpha_2}h_2$ satisfies $v_i(h'') \geq
v_i(h') \geq v_i(h)$ for all $i \in I$ and $v_j(h'') > v_j(h') =
v_j(h).$

If $h_2 = h$, then we must have $h'' = -c_1G^{\alpha_1}h_1$ and
$v_k(h'') = v_k(G^{\alpha_1}h_1) < v_k(h')$, which is a
contradiction. It follows that $h_2 \in H'$. So we have $v_i(h) \leq
v_i(h') \leq v_i(G^{\alpha_2}h_2)$ for all $i \in I$ and $v_j(h) =
v_j(G^{\alpha_2}h_2)$. Consequently, there exists $c'\in\mathbb{K}$ such that
$f'=f-c'G^{\alpha+\alpha_2}h_2$ is a $j$-reduction of $f$ modulo
$(H',G)$. Hence, $H'$ is a Standard Basis for $\mathcal{I}$. \cqd

Besides that as elements in a minimal Standard Basis have pairwise distinct values, we have the following result.

\begin{proposition} If $H$ and $H'$ are Standard Bases for $\mathcal{I}$ with
$H$ minimal, then $v(H)\subseteq v(H')$. In particular, all the
minimal Standard Bases for $\mathcal{I}$ have the same set of
values.
\end{proposition}

\pf Given $h \in H$ we will show that $v(h)\in v(H')$.

Since $H'$ is a Standard Basis for $\mathcal{I}$, if $k\in I_{h}$
there exist a $G$-product $G^\alpha$ and $h' \in H'$ such that
$v_i(h)\leq v_i(G^{\alpha}h')$ for all $i\in I$ and
$v_k(h)=v_k(G^{\alpha}h')$. In particular, $k\in I_{h'}$.

On the other hand, there exist a $G$-product $G^{\beta}$ and $h'' \in H$ such that $v_i(h')\leq v_i(G^{\beta}h'')$ for all $i\in I$ and $v_k(h')=v_k(G^{\beta}h'').$

In this way, we have $v_i(h) \leq v_i(G^{\alpha+\beta}h'')$, for all
$i \in I$ and $v_k(h) = v_k(G^{\alpha+\beta}h'')$. But since $H$ is
a minimal Standard Basis for $\mathcal{I}$, we must have $h''=h$ and
$\alpha =\beta = \underline{0}$.

Hence, $v_i(h) \leq v_i(h') \leq v_i(h)$ for all $i \in I$, that
is, $v(h)=v(h') \in v(H')$. \cqd

By the above proposition and Corollary \ref{finitgerad}, if $H$ is a
minimal Standard Basis for $\mathcal{I}$, then $v(H)$ is the unique
minimal system of generators for the $\Gamma$-semimodule of values
$v(\mathcal{I})$.

The following notion is the key for the algorithm that we will propose.

\begin{definition}
    Given $k \in I$, an $S_k$-\emph{process} of (the pair $(h_1,h_2)$ of) $H$ over $G$ is an element of the form
    $$c_1G^{\alpha_1}h_1 + c_2G^{\alpha_2}h_2,$$ where $c_1, c_2 \in \mathbb{K}$ and $G^{\alpha_1}$ and $G^{\alpha_2}$ are $G$-products in such way that $$v_k(c_1G^{\alpha_1}h_1 + c_2G^{\alpha_2}h_2) > \min \{v_k(G^{\alpha_1}h_1), v_k(G^{\alpha_2}h_2)\}.$$
\end{definition}

If $G=\{g_1,\ldots ,g_m\}$ then an $S_k$-process
$c_1G^{\alpha_1}h_1+c_2G^{\alpha_2}h_2$ of $h_1,h_2\in H$ over $G$
is determined, modulo an element in $\mathbb{K}\setminus\{0\}$, by a
solution $(\alpha_1,\alpha_2)\in\mathbb{N}^{2m}$ of the Diophantine
equation
$$\sum_{j=1}^{m}\alpha_{1j}v_k(g_j)+v_k(h_1)=\sum_{j=1}^{m}\alpha_{2j}v_k(g_j)+v_k(h_2).$$

The next result is a generalization of Theorem 4.1 of
\cite{coloquio} considering the product order $\leq$.

\begin{theorem}\label{teorema_equiv}
Let $H$ be a finite set of generators of $\mathcal{I}$. The following statements are
equivalent.

\noindent {\rm (a)} $H$ is a Standard Basis for $\mathcal{I}$.

\noindent {\rm (b)} Any nonzero $S_k$-process of $H$ over $G$ has a
vanishing final reduction modulo $(H,G)$.

\noindent {\rm (c)} Any nonzero $S_k$-process $c_1G^{\alpha_1}h_1 +
c_2G^{\alpha_2}h_2$ of $H$ over $G$ admits a representation as a sum
of the form $\sum_{\delta \in
\Delta}a_{\delta}G^{\delta}h_{\delta}$, where $\Delta \subset
\mathbb{N}^{\sharp G}$, $a_{\delta} \in \mathbb{K}$, $G^{\delta}$ is
a $G$-product, $h_{\delta} \in H$ and $$ \min\{v_i
(G^{\delta}h_{\delta}); \ \delta \in \Delta\} \geq
\min\{v_i(G^{\alpha_1}h_1), v_i(G^{\alpha_2}h_2)\},$$ for all $i \in
I$ and the inequality holds for $i =k$.

\end{theorem}

\pf (a) $\Rightarrow$ (b) It is immediate.

\noindent (b) $\Rightarrow$ (c) Let $f = c_1G^{\alpha_1}h_1 +
c_2G^{\alpha_2}h_2\neq 0$ be an $S_k$-process of $H$ over $G$ with a
vanishing final reduction modulo $(H,G)$, that is, $ f =
\sum_{\delta \in \Delta}a_{\delta}G^{\delta}h_{\delta}$, where
$\Delta\subset\mathbb{N}^{\sharp G}$, $a_{\delta} \in \mathbb{K}$,
$G^{\delta}$ is a $G$-product and $h_{\delta} \in H$. In particular, by the reduction process,
we have that $v(f) \leq v(G^{\delta}h_{\delta})$, for all $\delta
\in \Delta$ and $$ \min\{v_i (G^{\delta}h_{\delta}); \ \delta \in \Delta\} \geq v_i(f) \geq
\min\{v_i(G^{\alpha_1}h_1), v_i(G^{\alpha_2}h_2)\},$$ for all $i \in
I$. Furthermore, as $f$ is an $S_k$-process, we have  $$ \min\{v_k
(G^{\delta}h_{\delta}); \ \delta \in \Delta\} >
\min\{v_k(G^{\alpha_1}h_1), v_k(G^{\alpha_2}h_2)\}.$$

\noindent (c) $\Rightarrow$ (a) Let $f \in \mathcal{I}\setminus \{0\}$ and let $R_f$ be the set of all representations of $f$ as a sum of the form $f = \sum_{\delta \in \Delta}a_{\delta}G^{\delta}h_{\delta}$, where $\Delta \subset \mathbb{N}^{\sharp G}$, $a_{\delta} \in \mathbb{K}$, $G^{\delta}$ is a $G$-product and $h_{\delta} \in H$. Now consider the set $$\inf(R_f) : = \left\{ \inf\{ v(G^{\delta}h_{\delta}) ; \ \delta \in \Delta\} ; \ \ \sum_{\delta \in
    \Delta}a_{\delta}G^{\delta}h_{\delta} \in R_f \right\}.$$

This set is not empty and, since $\inf\{
v(G^{\delta}h_{\delta}); \ \delta \in \Delta\} \leq v(f)$,
$\inf(R_f)$ admits maximal elements. In this way, we can
consider a representation $ f = \sum_{\delta \in
\Delta}a_{\delta}G^{\delta}h_{\delta}$ such that $\lambda: =
\inf\{ v(G^{\delta}h_{\delta}); \ \delta \in \Delta\}$ is a
maximal element of $\inf(R_f)$.

If $\lambda = v(f)$, then by Remark \ref{Obs_Altura} and Remark
\ref{corolarioReducaoSB} we see that $f$ has a reduction modulo $(H, G)$.

Suppose by absurd that $\lambda < v(f)$. Setting $\lambda_i :
=\min\{ v_i(G^{\delta}h_{\delta}); \ \delta \in \Delta\}$, for $i\in
I$, we have $\lambda_k < v_k(f)$ for some $k \in I$. This means that
$\lambda_k = v_k(G^{\delta_1}h_{\delta_1}) = ... =
v_k(G^{\delta_s}h_{\delta_s})$ for some $\delta_1, ..., \delta_s \in
\Delta$ with $s \geq 2$ and $f' :=
a_{\delta_1}G^{\delta_1}h_{\delta_1} +
ba_{\delta_2}G^{\delta_2}h_{\delta_2}$ is an $S_k$-process of $H$
over $G$, for some $b\in\mathbb{K}$. Now, by hypothesis, there
exists a representation $f' = \sum_{\theta \in
\Theta}b_{\theta}G^{\theta}h_{\theta}$, where $\Theta \subset
\mathbb{N}^{\sharp G}$, $b_{\theta} \in \mathbb{K}$, $G^{\theta}$ is
a $G$-product, $h_{\theta} \in H$ and
\begin{equation}\label{desigualdade} \min\{v_i
(G^{\theta}h_{\theta}); \ \theta \in \Theta\} \geq
\min\{v_i(G^{\delta_1}h_{\delta_1}),
v_i(G^{\delta_2}h_{\delta_2})\},
\end{equation} for all $i \in I$ and the
inequality holds for $i =k$.

In this way, we have $a_{\delta_1}G^{\delta_1}h_{\delta_1} +
a_{\delta_2}G^{\delta_2}h_{\delta_2} = (1 -
b)a_{\delta_2}G^{\delta_2}h_{\delta_2} + \sum_{\theta \in
\Theta}b_{\theta}G^{\theta}h_{\theta}$. Then, we can write $$f =
(1-b)a_{\delta_2}G^{\delta_2}h_{\delta_2} + \sum_{\theta \in
\Theta}b_{\theta}G^{\theta}h_{\theta} + \sum_{\delta \in \Delta
\backslash \{\delta_1, \delta_2\}}a_{\delta}G^{\delta}h_{\delta}.$$

Set $\lambda' : = (\lambda'_1, ..., \lambda'_r)$, where $$\lambda_i' =  \min\{v_i(G^{\delta_2}h_{\delta_2}),
v_i(G^{\theta}h_{\theta}), v_i(G^{\delta}h_{\delta}); \ \theta \in
\Theta, \ \delta \in \Delta \setminus\{ \delta_1, \delta_2\}\}, \ \text{for} \ i\in I.$$

If $s \geq 3$, then $v_k(G^{\delta_s}h_{\delta_s}) = \lambda_k$ and
$\lambda'_k = \lambda_k$.

By (\ref{desigualdade}), we have $\lambda'_i \geq \lambda_i$ for all $i\in I\setminus\{k\}$. If the inequality holds for some $i \in
I\backslash \{k\}$, we obtain $\lambda' > \lambda,$ which is an
absurd because $\lambda$ is a maximal element of $\inf(R_f)$.

On the other hand, if $\lambda'_i = \lambda_i$, for all $i \in
I\setminus\{k\}$, we have $\lambda' = \lambda$ with  $\lambda'_k =
v_k(G^{\delta_j}h_{\delta_j}) = ... =
v_k(G^{\delta_{s}}h_{\delta_{s}})$, for $j \geq 2$.

In this way, repeating the above argument, we can suppose that $s = 2$.

If $b \neq 1$, then
$\lambda_k<v_k(f)=\lambda'_k = v_k(G^{\delta_2}h_{\delta_2}) =
\lambda_k$, which is a contradiction.

If $b = 1$, we have $\lambda'_k > \lambda_k$ and $\lambda'_i = \min\{v_i(G^{\theta}h_{\theta}),
v_i(G^{\delta}h_{\delta}); \ \theta \in \Theta, \ \delta \in \Delta \setminus
\{\alpha, \beta\}\}$ for $i \in
I\backslash\{k\}$. In this way, we obtain $\lambda'_i \geq
\lambda_i$, for all $i \in I\backslash\{k\}$. Hence, $\lambda' >
\lambda$, but this is an absurd because $\lambda$ is a maximal
element of $\inf(R_f)$.

Therefore, we must have $\lambda = v(f)$ and $H$ is a Standard basis
for $\mathcal{I}$. \cqd

The characterization given by the item (b) of Theorem
\ref{teorema_equiv} allows us to obtain an algorithm to compute a
Standard Basis for $\mathcal{I}$.

\begin{theorem} \label{teoremaAlgorimo}
Let $H_0$ be a finite set of generators for $\mathcal{I}$ such that
$\cup_{i\in I}B_i \subseteq H_0$, where $B_i$ is described in
Theorem \ref{existBase2}. Then we always obtain a Standard Basis $H$ for $\mathcal{I}$ with the following
algorithm:
\end{theorem}

\begin{center}
\textbf{ALGORITHM 1.} Standard Basis for $\mathcal{I}$
\vspace{0.3cm}

\begin{tabular}{|l|}
\hline
\textbf{input}: $G, H_0$; \\

\textbf{define}: $H_{-1} := \emptyset$ \textbf{and}    $j := 0$; \\

\textbf{while} $H_{j} \neq H_{j - 1}$ \textbf{do} \\

\hspace{9mm} $\rho:=$ \text{an upper bound for} $\kappa=(\kappa_1,\ldots ,\kappa_r)$; \\

\hspace{9mm} $\mathcal{S} := \cup_{k\in I}\{f; \ f \ \text{is an} \
S_k\text{-process of} \ H_{j} \ \text{over} \ G \ \text{and} \
v_i(f) < \rho_i \ \text{for some} \ i \in I\}$; \\

\hspace{9mm} $\mathcal{R} := \{h; \ h\neq 0 \ \text{is a final reduction
of} \ f \in \mathcal{S}\ \text{modulo} \ (H_j,G) \}$; \\

\hspace{9mm} $H_{j+1} := H_{j} \cup \mathcal{R}$; \\

\textbf{output}: $H = H_{j+1}$. \\

\hline
\end{tabular}
\end{center}

\vspace{2mm}

\pf Notice initially that in each iteration the set $\mathcal{S}$ is
finite\footnote{Here we identify two $S_k$-processes $f$ and $f_1$
if $f=cf_1$ for some $c\in\mathbb{K}\setminus\{0\}$.} due to
the condition ``$v_i(f) < \rho_i$ for some $i \in I$" imposed on an
$S_k$-process $f$ of $H_j$ over $G$.

In fact, for each pair of elements $(h,h')$ of $H_j$ let us consider
$S_k$-processes on the form $f = cG^{\alpha}h + c'G^{\alpha'}h'$ of
$H_j$ over $G= \{g_1, ..., g_m\}$ such that $v_i(f) =: \lambda_i<
\rho_i$ for some $i \in I$. Without loss of generality we can suppose that $i=1$.

If $f$ is not an $S_1$-process, then $\lambda_1 =
\min\{v_1(G^{\alpha}h),v_1(G^{\alpha'}h')\}$ and we must have
$$\sum_{l=1}^{m}\alpha_l\cdot v_1(g_l) + v_1(h) = \lambda_1\
\ \ \mbox{or}\ \ \ \sum_{l=1}^{m}\alpha'_l\cdot v_1(g_l) + v_1(h') =
\lambda_1.$$ As such Diophantine equations have a finite number of non-negative solutions we obtain a finite number of $S_k$-processes $f$
with $v_1(f)=\lambda_1$.

On the other hand, if $f$ is also an $S_1$-process, then we have
that $v_1(G^{\alpha}h)=v_1(G^{\alpha'}h')=s<\lambda_1$.
Consequently,
$$\sum_{l=1}^{m}\alpha_l\cdot v_1(g_l) + v_1(h) = s=\sum_{l=1}^{m}\alpha'_l\cdot v_1(g_l) + v_1(h')$$ and we find a finite number of
$S_k$-processes $f$ with $v_1(f)=\lambda_1$.

The hypothesis ``$\cup_{i\in I}B_i \subset H_0$'' guarantees
that if $v(f) \geq \rho\geq\kappa$, then $f$ admits a vanishing
final reduction modulo $(H_0,G)$ and, consequently, modulo $(H_j,
G)$ for all $j \in \mathbb{N}$.

Now consider $H = \cup_{j \geq 0}H_j$ and let $f$ be an
$S_k$-process of $H$ over $G$. In this way, $f$ is also an
$S_k$-process of $H_j$ over $G$ for some $j \geq 0$. By the
algorithm above, $f$ has a vanishing final reduction modulo
$(H_{j+1},G)$. Hence, its final reduction modulo $(H,G)$ is
also zero.

Finally we will show that $H$ is finite.

If $h \in H\setminus H_0$, then $v_i(h) < \rho_i$ for all $i \in I_h$ because otherwise, as $H_0\supseteq \cup_{i\in I}B_i$, $h$ would have a reduction module $(H_0,G)$. Therefore there exist finite possibilities for $v(h)$, that is, $v(H)$ is finite. Now suppose, by absurd, that $H$ is not finite, then for all $j \in \mathbb{N}$ there exists at least an element $h_j \in H$ such that $h_j \in H_{j+1}\setminus H_j$. In this way, we must have $v(h_j) \not\in v(H_j)$, for all $j$, which makes the set $v(H)$ infinite and, consequently, a contradiction.    \cqd

In order to obtain a Standard Basis for $\mathcal{O}$, it is
sufficient to make a few modifications in Theorem
\ref{teoremaAlgorimo}, in the previous algorithm and in its proof.
More precisely, we consider a finite set of generators $G_0$ for
$\mathcal{O}$ (as $\mathbb{K}$-algebra) containing the sets $B_i$
for $i\in I$ described in Remark \ref{Conjuntos B_i's}.

\begin{center}
    \textbf{ALGORITHM 2.} Standard Basis for $\mathcal{O}$
    \vspace{0.3cm}

    \begin{tabular}{|l|}
        \hline
        \textbf{input}: $G_0$; \\

        \textbf{define}: $G_{-1} := \emptyset$ \textbf{and}    $j := 0$; \\

        \textbf{while} $G_{j} \neq G_{j - 1}$ \textbf{do} \\

        \hspace{9mm} $\mathcal{S} := \cup_{k\in I}\{g; \ g \ \text{is an} \
        S_k\text{-process of} \ G_{j} \ \text{and} \
        v_i(g) < \sigma_i \ \text{for some} \ i \in I\}$; \\

        \hspace{9mm} $\mathcal{R} := \{h; \ h\neq 0 \ \text{is a final reduction
            of} \ g \in \mathcal{S}\ \text{modulo} \ G_j
        \}$; \\

        \hspace{9mm} $G_{j+1} := G_{j} \cup \mathcal{R}$; \\

        \textbf{output}: $G = G_{j+1}$. \\

        \hline
    \end{tabular}
\end{center}

\vspace{2mm}

Notice that, in Algorithm 1 (and consequently in Algorithm 2), if we consider only $S_k$-processes and
$k$-reductions for a fixed $k\in I$ we obtain a Standard Basis
$H_{\{k\}}$ for $\mathcal{I}_k$. Moreover, every $S_k$-process of
$h_1,h_2\in H_{\{k\}}$ over $G$ has a vanishing $k$-reduction module
$(H_{\{k\}},G)$, that is, $v_k\left (
c_1G^{\alpha_1}h_1+c_2G^{\alpha_2}h_2-\sum_{i=1}^{n_k}g_ih_i\right
)=\infty$ where $g_i$ is a $\mathbb{K}$-combination of $G$-products
and $n_k=\sharp H_{\{k\}}$.

In this way
\begin{equation}\label{syzy}(g_1-c_1G^{\alpha_1},g_2-c_2G^{\alpha_2},g_3,\ldots
,g_{n_k})\end{equation} is a syzygy of $(h_1,\ldots ,h_{n_k})$
considering the $\mathcal{O}_k$-module structure. By the Schreyer
Theorem (see Theorem 15.10 in \cite{einsenbud}), the first module of
syzygies of $(h_1,\ldots ,h_{n_k})$ is generated by the elements as
in (\ref{syzy}) and then we have a set of generators for the
$\mathcal{O}$-module
$\{f\in\mathcal{I};\ v_k(f)=\infty\}$.

By successive computations as described, we can obtain a set of
generators for $\mathcal{T}^{i}$ for $i\in I$ and a set $B_i$ as
mentioned in the hypothesis of Theorem \ref{teoremaAlgorimo}.

Hence, in theory, we can obtain a
Standard Basis for a fractional ideal $\mathcal{I}$ of $\mathcal{O}$ applying the Algorithm 1 for a finite set of generators $H_0$ of $\mathcal{I}$.

We illustrate the above remark with the following example.

\begin{example}

Consider the space curve $Q = P_1 \cap P_2 \subset
\mathbb{C}[[X,Y,Z]]$, where the prime ideal $P_1$ is generated by
$$f_1=3Z^2 - 4X^2Z+XY^2-3X^2Y-4X^4\ \ \ \mbox{and}$$
$$f_2=24Y^2Z-18XYZ-224X^3Z+9Y^3+32X^2Y^2-96X^3Y-128X^5+9X^4$$ and
$P_2$ is generated by
$$g_1= YZ^2+XY^3-2X^2YZ-2X^2Y^2+X^3Z-X^4Y,\ \ \ g_2=2XZ^2-Y^2Z+X^3Z-X^3Y$$
$$g_3=Z^3-XY^3\ \ \mbox{and}\ \ g_4=2Y^3+XZ^2-3XYZ+X^2Y^2-2X^3Z-4X^3Y+X^4-X^5.$$

Notice initially that $\varphi_1(t_1) = (t_1^6, t_1^8+2t_1^9, t_1^{10}+t_1^{11})$ and $\varphi_2(t_2) = (t_2^6, t_2^8+t_2^9, t_2^{10}+t_2^{11})$ are parameterizations of the branches $P_1$ and $P_2$ respectively. 

We consider $G_0=\{x,y,z,f_1,f_2,g_1,g_2,g_3,g_4\}$ in the Algorithm 2. As we remarked, it can be taken just $\{x,y,z\}$ but in this case we will have several steps in the algorithm.

Computing $S_k$-processes of $G_0$ we find $\{h_1=y^2-xz,h_2=yz-x^3,h_3=z^2-x^2y,g_5=xg_1+g_4,g_6=yg_1-xg_2,f_3=3zf_1+xf_2,f_4=3x^2f_1+yf_2,f_5=3xyf_1+zf_2\}$.
All these elements coincide with their final reduction modulo $G_0$. In the next step we consider $G_1=G_0\cup\{h_1,h_2,h_3,g_5,g_6,f_3,f_4,f_5\}$. Since all $S_k$-process of $G_1$ has a vanishing final reduction, the algorithm stops and $G_1$ is a Standard Basis for $\mathcal{O}$.

Notice that $v(g_4)=v(xg_1)$ and $v(h_3)=v(x^2z+f_1)$, so we can discard $g_4$ and $h_3$ from $G_1$ in such way that we obtain, by Proposition \ref{minima}, a minimal Standard Basis
$$G = \{x,y,z,f_1,f_2,f_3,f_4,f_5,g_1,g_2,g_3,g_5,g_6,h_1,h_2\},$$ its respective set of values
\begin{center} $v(G)=\{(6,6),(8,8),(10,10),(\infty,21),(\infty,25),(\infty,32),(\infty,34),(\infty,36)$,

\hspace{2cm} $(25,\infty),(27,\infty),(29,\infty),(31,\infty),(32,\infty),(34,\infty),(17,17),(19,19)\}$
\end{center}
is the minimal set of generators of the semiring $\Gamma$ and its conductor is $(31,31)$.

Notice that $(25,26)\in\Gamma$. In fact we have
$$(25,26)=v(g_3+xz^2)=\inf\{(25,\infty),(6,6)+2(10,10)\}=(25,\infty)\oplus \left((6,6)\odot (10,10)^2\right).$$

By the other hand, using (\ref{elemento}), we can verify that $(24,26)\not\in\Gamma$.

As we remarked Standard Bases for $\mathcal{O}_i$ and $Q^i$ are obtained in the steps of the algorithm with input a set of generators for $\mathcal{O}$. In fact, if we apply the algorithms presented in \cite{Abramo-Marcelo}, we obtain that
$B'_1 = \{x,y,z,h_1,h_2\}$ and $B'_2=\{x,y,z,h_1,h_2,h_3\}$ are Standard Bases for $\mathcal{O}_1$ and $\mathcal{O}_2$ and its respective semigroups are $\Gamma_1\cap\mathbb{N}=\langle 6,8,10,17,19\rangle$ and $\Gamma_2\cap\mathbb{N}=\langle 6,8,10,17,19,21\rangle$. In addition,
$B''_1 = \{g_1,g_2,g_3,g_4,g_5,g_6\}$ and $B''_2=\{f_1,f_2,f_3,f_4,f_5\}$
are Standard Bases for $Q^1$ and $Q^2$ respectively. We have $v_1(Q^1)=\{25,27,29,31,32,34\}+\Gamma_1$, $v_2(Q^2)=\{21,25,32,34,36\}+\Gamma_2$, the conductors of $v_1(Q^1)$ and $v_2(Q^2)$ are $\sigma_1=\sigma_2=31$. 

In the next diagram we indicate the common elements in $\Gamma$ and the box $[0,31]\times [0,31]$.
\end{example}
\vspace{0.5cm}

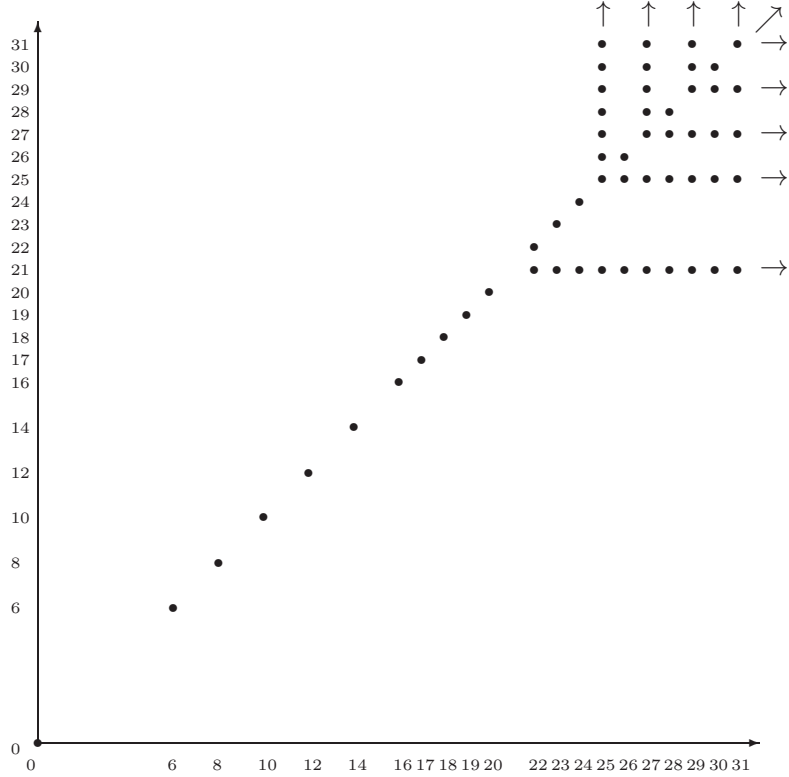
\begin{center}
\setlength{\unitlength}{0.3cm}
\begin{picture}(32,32)
\put(0,0){\vector(1,0){32}}
\put(0,0){\vector(0,1){32}}
\put(0,0){\circle*{0.3}}
\put(6,6){\circle*{0.3}}
\put(8,8){\circle*{0.3}}
\put(10,10){\circle*{0.3}}
\put(12,12){\circle*{0.3}}
\put(14,14){\circle*{0.3}}
\put(16,16){\circle*{0.3}}
\put(17,17){\circle*{0.3}}
\put(18,18){\circle*{0.3}}
\put(19,19){\circle*{0.3}}
\put(20,20){\circle*{0.3}}
\put(22,22){\circle*{0.3}}
\put(23,23){\circle*{0.3}}
\put(24,24){\circle*{0.3}}
\put(25,25){\circle*{0.3}}
\put(26,26){\circle*{0.3}}
\put(27,27){\circle*{0.3}}
\put(28,28){\circle*{0.3}}
\put(29,29){\circle*{0.3}}
\put(30,30){\circle*{0.3}}
\put(31,31){\circle*{0.3}}
\put(32,30.75){$\rightarrow$}
\put(31.75,31.75){$\nearrow$}
\put(30.75,32){$\uparrow$}
\put(32,20.75){$\rightarrow$}
\put(22,21){\circle*{0.3}}
\put(23,21){\circle*{0.3}}
\put(24,21){\circle*{0.3}}
\put(25,21){\circle*{0.3}}
\put(26,21){\circle*{0.3}}
\put(27,21){\circle*{0.3}}
\put(28,21){\circle*{0.3}}
\put(29,21){\circle*{0.3}}
\put(30,21){\circle*{0.3}}
\put(31,21){\circle*{0.3}}
\put(32,24.75){$\rightarrow$}
\put(26,25){\circle*{0.3}}
\put(27,25){\circle*{0.3}}
\put(28,25){\circle*{0.3}}
\put(29,25){\circle*{0.3}}
\put(30,25){\circle*{0.3}}
\put(31,25){\circle*{0.3}}
\put(32,26.75){$\rightarrow$}
\put(28,27){\circle*{0.3}}
\put(29,27){\circle*{0.3}}
\put(30,27){\circle*{0.3}}
\put(31,27){\circle*{0.3}}
\put(32,28.75){$\rightarrow$}
\put(30,29){\circle*{0.3}}
\put(31,29){\circle*{0.3}}
\put(24.75,32){$\uparrow$}
\put(25,26){\circle*{0.3}}
\put(25,27){\circle*{0.3}}
\put(25,28){\circle*{0.3}}
\put(25,29){\circle*{0.3}}
\put(25,30){\circle*{0.3}}
\put(25,31){\circle*{0.3}}
\put(26.75,32){$\uparrow$}
\put(27,28){\circle*{0.3}}
\put(27,29){\circle*{0.3}}
\put(27,30){\circle*{0.3}}
\put(27,31){\circle*{0.3}}
\put(28.75,32){$\uparrow$}
\put(29,30){\circle*{0.3}}
\put(29,31){\circle*{0.3}}
{\tiny \put(-0.5,-1.2){0}\put(5.75,-1.2){6}\put(7.75,-1.2){8}\put(9.75,-1.2){10}\put(11.75,-1.2){12}\put(13.75,-1.2){14}\put(15.75,-1.2){16}\put(16.75,-1.2){17}\put(17.75,-1.2){18}
\put(18.75,-1.2){19}\put(19.75,-1.2){20}\put(21.75,-1.2){22}\put(22.75,-1.2){23}\put(23.75,-1.2){24}\put(24.75,-1.2){25}
\put(25.75,-1.2){26}\put(26.75,-1.2){27}\put(27.75,-1.2){28}\put(28.75,-1.2){29}\put(29.75,-1.2){30}\put(30.75,-1.2){31}
}
{\tiny \put(-1.2,-0.5){0}\put(-1.2,5.75){6}\put(-1.2,7.75){8}\put(-1.2,9.75){10}\put(-1.2,11.75){12}\put(-1.2,13.75){14}\put(-1.2,15.75){16}\put(-1.2,16.75){17}\put(-1.2,17.75){18}
\put(-1.2,18.75){19}\put(-1.2,19.75){20}\put(-1.2,20.75){21}\put(-1.2,21.75){22}\put(-1.2,22.75){23}\put(-1.2,23.75){24}\put(-1.2,24.75){25}
\put(-1.2,25.75){26}\put(-1.2,26.75){27}\put(-1.2,27.75){28}\put(-1.2,28.75){29}\put(-1.2,29.75){30}\put(-1.2,30.75){31}
}
\end{picture}
\end{center}
\vspace{0.5cm}

\section{The Module of K\"{a}hler Differentials for Plane Curves}

In this section, we will consider the module of K\"{a}hler
differentials $\Omega_{\mathcal{O}/\mathbb{K}}$ for an algebroid
plane curve. In the analytical case,
$\Omega_{\mathcal{O}/\mathbb{C}}$ is an important example of
fractional ideal since the relative ideal associated to it plays a
central role in the analytic classification problem, as we can see
in \cite{Abramo-Marcelo3} for irreducible plane curves and in
\cite{AME} for plane curves with two branches.

In the sequel, $\mathbb{K}$ will be an algebraically closed field of characteristic zero and $Q = \langle f \rangle$ denotes an algebroid reduced plane curve, with $f = \prod_{i \in I}f_i \in \mathbb{K}[[X,Y]]$, $\langle f_i \rangle \neq \langle f_j \rangle$ and $f_i$ irreducible for any $i \in I=\{1,\ldots ,r\}$.

\begin{definition}
    The \emph{module of K\"{a}hler differentials} over $\mathcal{O}$
    is the $\mathcal{O}$-module $$\Omega:=\Omega_{\mathcal{O}/\mathbb{K}} = \frac{\mathcal{O}dx+\mathcal{O}dy}{(f_xdx+f_ydy)\mathcal{O}}.$$
\end{definition}

In the same way, for every $i \in I$, the module of K\"{a}hler
differentials over $\mathcal{O}_i$ is $$\Omega_i :=
\frac{\mathcal{O}_idx+\mathcal{O}_idy}{((f_i)_xdx+(f_i)_ydy)\mathcal{O}_i}.$$

Given $\omega_i = p_{i}dx + q_{i}dy \in \Omega_i$, we define
$$\varphi_i^* (p_{i}dx + q_{i}dy) =
t_i\cdot\left( p_{i}(\varphi_i(t_i))\cdot x'(t_i) +
q_{i}(\varphi_i(t_i))\cdot y'(t_i)\right ) \in
\overline{\mathcal{O}_i} = \mathbb{K}[[t_i]],$$ where
$\varphi_i(t_i) = (x(t_i), y(t_i))$ is a parameterization of the
branch $\langle f_i \rangle$ and we denote $\varphi_i^*(\Omega_i) =
\{\varphi_i^*(\omega_i); \ \omega_i \in \Omega_i\}$.

Furthermore, we consider the $\mathcal{O}$-modules homomorphism
$\varphi^* : \Omega \rightarrow \overline{\mathcal{O}}$ defined by \begin{align} & \Omega \ \longrightarrow \ \ \ \ \bigoplus_{i\in I}\Omega_i \ \ \ \hspace{0.5mm} \longrightarrow \ \ \ \ \overline{\mathcal{O}} \cong \bigoplus_{i\in I}\mathbb{K}[[t_i]] \n \\ & \omega \ \longmapsto \ (\omega_1, ..., \omega_r) \ \longmapsto \ (\varphi_1^*(\omega_1), ..., \varphi_r^*(\omega_r))\n \end{align} where $\omega = pdx + qdy \in \Omega$ and $\omega_i = p_{i}dx + q_{i}dy \in \Omega_i$, for all $i \in I.$

In what follows we will denote by $\varphi^*(\Omega)$ the image of $\Omega$ by the previous $\mathcal{O}$-module homomorphism $\varphi^*$. Remark that $\{dh:= h_xdx + h_ydy; \ h \in \mathcal{O}\} \subset \Omega$ and if $g \in \mathcal{O}$ is such that $v(g) = \sigma$, we get $g\varphi^*(\Omega) \subset \mathcal{O}$. Therefore, $\varphi^*(\Omega)$ can be considered as a fractional ideal of $\mathcal{O}$.

For each $i \in I$, we define the value of an element $\omega_i \in \Omega_i$ as $\nu_i(\omega_i) = v_i(\varphi_i^*(\omega_i)),$ where $v_i$ is the discrete normalized valuation of $\overline{\mathcal{O}}_i$, and the set of the values of the differentials of $\mathcal{O}_i$ as $$\Lambda_i := \left\{\nu_i(\omega_i); \ \omega_i \in \Omega_i \right\} \subset \overline{\mathbb{N}}.$$

We naturally define the value $\nu(\omega): = (\nu_1(\omega_1), ...,
\nu_r(\omega_r))$ of the element $\omega \in \Omega$ and we write
$$\Lambda = \left\{\nu(\omega); \ \omega \in \Omega\right\} \subset
\bigoplus_{i\in I}\Lambda_i \subset \overline{\mathbb{N}}^r$$ for
representing the relative ideal associated to $\Omega$. Notice that
the set $\Lambda$ can be obtained by a Standard Basis for
$\varphi^*(\Omega)\subset\overline{\mathcal{O}}$ that will be called a Standard Basis for $\Omega$.

We remark that in \cite{Abramo-Marcelo3} and \cite{AME} the authors considered  $\Lambda\cap \mathbb{N}^r$ $(r=1,2)$ as the main ingredient to proceed an answer to the analytic classification problem for curves with one and two branches. In addition, as we mentioned in Introduction, the set $\Lambda$ is related to the values of the module of logarithmic residues along a complete intersection curve $Q$ and to the set of values of the Jacobian ideal of $Q$.

For the irreducible case (plane or not), the set $\Lambda$ can be computed using the algorithms in \cite{Abramo-Marcelo}. For plane curves with two branches, Pol presents in Subsection 4.3.3 of \cite{pol} a method that allows to compute $\Lambda$.

The next result is a generalization of Algorithm 4.10 in
\cite{Abramo-Marcelo} in the sense that we can compute a Standard
Basis for $\Omega$ by a Standard Basis $G$ for the local ring
$\mathcal{O}$ such that $\{f_i;\ i\in I\}\subset G$.

\begin{proposition} \label{teoremaAlgorimo3}
    Let $G$ be a Standard Basis for $\mathcal{O}$ with $\{f_i;\ i\in I\}\subseteq G$. We
    always obtain a Standard Basis $H$ for the module of K\"{a}hler differentials $\Omega$
with the following algorithm:
\end{proposition}

\begin{center}
    \textbf{ALGORITHM 3.} Standard Basis for $\Omega$
    \vspace{0.3cm}

    \begin{tabular}{|l|}
        \hline
        \textbf{input}: $G$; \\

        \textbf{define}: $H_{-1} := \emptyset$, $H_0 := \{dg; \ g \in G \}$ \textbf{and}    $j := 0$; \\

        \textbf{while} $H_{j} \neq H_{j - 1}$ \textbf{do} \\

        \hspace{9mm} $\mathcal{S} := \cup_{k=1}^{r}\{\omega; \ \omega \ \text{is an} \
        S_k\text{-process of} \ H_{j} \ \text{over} \ G \ \text{and} \
        \nu_i(\omega) < \sigma_i \ \text{for some} \ i \in I\}$; \\

        \hspace{9mm} $\mathcal{R} := \{\varpi; \ \varpi\neq 0 \ \text{is a final reduction
            of} \ \omega\in \mathcal{S} \ \text{modulo} \ (H_j,G)
         \}$; \\

        \hspace{9mm} $H_{j+1} := H_{j} \cup \mathcal{R}$; \\

        \textbf{output}: $H = H_{j+1}$. \\

        \hline
    \end{tabular}
\end{center}

\vspace{2mm}

\pf Notice that the main difference between Algorithm 1 and
Algorithm 3 is the set of generators $H_0$ for $\Omega$. In
Algorithm 1 the hypothesis $\cup_{i\in I}B_i\subset H_0$ guarantees
that if $\nu(\omega) \geq \kappa$, then $\omega$ admits a vanishing
final reduction modulo $(H_0,G)$ and, consequently, modulo $(H_j,
G)$ for all $j \in \mathbb{N}$. We will show that this claim remains
true if we change $\kappa$ by $\sigma$.

In order to do this, it is sufficient to show that if $\nu_k(\omega)
\geq \sigma_k$, for some $k \in I_{\omega}$, then $\omega$ admits a
$k$-reduction modulo $(H_0,G)$, where $H_0=\{dg; \ g \in G \}$. For
simplicity, we will suppose $k=1$.

As $\sigma_1 = \mu_1 + \sum_{l=2}^{r}v_1(f_l)$, where $\mu_1$ is the conductor of $\Gamma_1$, the condition $\sigma_1 \leq \nu_1(\omega) < \infty $ implies that there exists a $G$-product $G^{\alpha}$ such that $\nu_1(\omega) = v_1(G^{\alpha}) + \sum_{l=2}^{r}v_1(f_l)$.

If $\alpha \neq \underline{0}$ then $v_1(G^{\alpha})=\nu_1(G^{\beta}dg)$ for some $g \in G$ and $$\nu_1(\omega) = \nu_1\left(\prod_{l=2}^{r}f_lG^{\beta}dg\right) = \nu_1(G^{\theta}dg).$$

If $\alpha =\underline{0}$ then $\mu_1=0$ and $\nu_1(\omega) =\nu_1\left(\prod_{l=2}^{r}f_l\right)$. So, $\nu_1(\omega) =\nu_1(G^{\theta}dg)$, where $G^{\theta}=\prod_{l=3}^{r}f_l$ and $dg=df_2$.

In any case $\nu_i(\omega) \leq \infty = \nu_i(G^{\theta}
dg)$ for $i =2,...,r$ and consequently, $\omega$ admits a $1$-reduction modulo
$(H_0,G)$.

As $\nu(\omega-c G^{\theta} dg)>\nu(\omega)\geq \sigma$ for some
$c\in\mathbb{C}$, we can repeat the same argument for $\omega-c
G^{\theta} dg$. Hence, $\omega$ admits a vanishing final reduction
modulo $(H_0,G)$.  \cqd

In the following example we apply Algorithm 3 in order to compute the minimal set of generators of $\Lambda$ for a plane curve with three branches.

\begin{example}

Consider the plane curve $Q = \langle
Y(Y-X^n)(Y-X^{m}-aX^{m+1})\rangle \subset \mathbb{C}[[X,Y]]$, with
$1<n <m$ and any $a\in \mathbb{C}$. Writing $f_1 = Y$, $f_2=Y-X^n$
and $f_3 = Y-X^{m}-aX^{m+1}$, then $\varphi_1(t_1)=(t_1,0)$,
$\varphi_2(t_2)=(t_2,t_2^{n})$ and
$\varphi_3(t_3)=(t_3,t_3^{m}+at_3^{m+1})$ are parameterizations for
the branches $\langle f_1 \rangle$, $\langle f_2 \rangle$ and
$\langle f_3 \rangle$ respectively.

It is not difficult to see (by Algorithm 2) that $G=\{x,y,f_2,f_3\}$
is a Standard Basis for $\mathcal{O}$, the semiring $\Gamma$ is
minimally generated by $$\{(1,1,1),(\infty ,n,m),(n,\infty ,
n),(m,n,\infty)\}$$ and the conductor of $\Gamma$ is $\sigma
=(n+m,2n,n+m)$.

Applying Algorithm 3 with $H_0:=\{dx, dy, df_2, df_3\}$ we perform
several computations involving $S_k$-processes and $k$-reductions,
some of them with infinite many steps. For this reason we omit the
iterations of the algorithm. As result of the computations we obtain that the set $H =\{dx, dy, df_2, df_3, \omega_1, \omega_2, \omega_3 \},$ where
$$\omega_1 = -nydx+xdy,$$
$$\omega_2=(m-n)xdy-m\omega_1+\sum_{i=1}^{\infty}\frac{(-1)^{i-1}a^in(m-n+1)^{i-1}}{(m-n)^i}x^i\omega_1,$$
$$\omega_3=(n-m)xdf_3+\omega_2-\frac{m}{n}x^{m-n}\omega_2-\frac{a(m+1)}{n}x^{m-n+1}\omega_2,$$
is such way that $H$ is a minimal Standard Basis for the module of
K\"{a}hler differentials $\Omega$ over the local ring of the curve
$Q$, the set
$$v(H)=\{(1,1,1), (\infty, n, m), (n, \infty, n), (m, n, \infty),
(\infty, \infty, m+1), (\infty, n+1, \infty), (m+1, \infty, \infty)
\}$$ is the minimal set of generators of the set $\Lambda$ and its
conductor is $\varrho=(m+1,n+1,m+1)$.

\end{example}

\section*{Acknowledgment}

The authors are grateful to the anonymous referee for the suggestions that improve this work.

\vspace{1cm}

\begin{tabular}{lcl}
    Carvalho, E. & & Hernandes, M. E. \\
    {\it emilio.carvalho$@$gmail.com} & & {\it mehernandes$@$uem.br} \\
    & DMA-UEM &  \\
    & Av. Colombo 5790 & \\
    & Maring\'{a}-PR 87020-900 & \\
    & Brazil &
\end{tabular}

\end{document}